# Derivation of new Degrees for Best (Co)convex and Unconstrained Polynomial Approximation in $\mathbb{L}_p^{\alpha,\beta}$ space: I

Malik Saad Al-Muhja        Habibulla Akhadkulov        Nazihah Ahmad

ABSTRACT.

The purpose of this work is to present the derivation and an estimate of the degrees of the best approximation based on convex, coconvex and unconstrained polynomials, and discuss some applications. We simplify the term convex and coconvex polynomial as (co)convex polynomial herein.

## 1. Introduction and Preliminaries

The (co)convex and unconstrained polynomial approximation (COCUNP approximation) is one of the important estimates in approximation theory with Jacobi weights that Kopotun has recently introduced (see [10, 11, 12, 15, 16, 19]), and that are defined as

$$w_{\alpha,\beta}(x) = (1+x)^\alpha (1-x)^\beta,$$

and

$$\mathbb{L}_p^{\alpha,\beta} = \left\{ f:[-1,1] \to \mathbb{R} : \left\| w_{\alpha,\beta} f \right\|_p < \infty, \text{ and } 0 < p < \infty \right\}.$$

$\mathbb{L}_p^{\alpha,\beta}$ space expansions will be used in the derivation of the degrees of the best unconstrained polynomial (UNP) approximation. Then the resulting degree after substitution will be rewriting the polynomial to (co)convex polynomial (COCP). Then, find the degree of the best coconvex polynomial (CCP) approximation which indicates two of inflection points at least and with one inflection point by using $\mathbb{L}_p^{\alpha,\beta}$ space, which will be explained in this paper. In fact, we also need the nomenclatures: unconstrained and convex polynomial (UNCP) approximation and Ditzian-Totik modulus of smoothness (DTMS). Next, let $\Delta_h^k(f,x)$ be the $k$th symmetric difference of $f$ is given [7] by

$$\Delta_h^k(f,x) = \begin{cases} \sum_{i=0}^k \binom{k}{i}(-1)^{k-i} f\left(x + \left(\frac{2i-k}{2}\right)h\right), & x \pm \frac{kh}{2} \in [-1,1] \\ 0, & \text{otherwise.} \end{cases}$$

**Definition 1.1 [17]** The space $L_p([-1,1])$, $0 < p < \infty$, denotes the space of all measurable functions $f$ on $[-1,1]$, such that

$$\|f\|_{L_p[-1,1]} = \begin{cases} \left(\int_{-1}^{1} |f(x)|^p dx\right)^{\frac{1}{p}} < \infty, & 0 < p < \infty \\ \underset{x \in [-1,1]}{\mathrm{esssup}} |f(x)|, & p = \infty. \end{cases}$$

Let $\|\cdot\|_p = \|\cdot\|_{L_p[-1,1]}$, $0 < p \le \infty$ and $\phi(x) = \sqrt{1-x^2}$. Then, DTMS of a function $f \in L_p[-1,1]$, is defined [6] by

$$\omega_{k,r}^{\phi}(f,t)_p = \sup_{0 < h \le t} \|\phi^r \Delta_{h\phi}^k(f,x)\|_p, \qquad k, r \in \mathbb{N}_o.$$

Also, the $k$th usually modulus of smoothness of $f \in L_p[-1,1]$ is defined [7] by

$$\omega_k(f,\delta,[-1,1])_p = \sup_{0 < h \le \delta} \|\Delta_h^k(f,x)\|_p, \qquad \delta > 0, p \le \infty.$$

Denote by $AC_{loc}(-1,1)$ and $AC[-1,1]$ the set of functions whose are locally absolutely continuous on $(-1,1)$ and absolutely continuous on $[-1,1]$ respectively. Now, we will need to accept the following:

**Definition 1.2 [20]** Let

$$\alpha, \beta \in J_p = \begin{cases} \left(\frac{-1}{p}, \infty\right), & \text{if } p < \infty, \\ [0, \infty), & \text{if } p = \infty. \end{cases}$$

Define

$$\mathbb{L}_{p,r}^{\alpha,\beta} = \left\{ f: [-1,1] \to \mathbb{R} : f^{(r-1)} \in AC_{loc}(-1,1), 1 \le p \le \infty \text{ and } \|w_{\alpha,\beta} f^{(r)}\|_p < \infty \right\}$$

and for convenience denote $\mathbb{L}_{p,0}^{\alpha,\beta} = \mathbb{L}_p^{\alpha,\beta}$.

**Definition 1.3 [24]** Let $f \in \mathbb{C}$ be a continuous function on $[-1,1]$, $r \ge 1$ and $\phi(x) = \sqrt{1-x^2}$. Denote by $\mathbb{C}^{(r)}$ the space of continuous functions which possess an absolutely continuous $(r-1)$st derivative in $[-1,1]$ such that $f^{(r)}$ is almost everywhere bounded; that is,

$$\mathbb{C}^{(r)} = \{f \in \mathbb{C} : f^{(r-1)} \in AC[-1,1] \text{ and } \|f^{(r)}\| < \infty\}.$$

Denote by $\ddot{\mathbb{C}}^{(r)}$ the space of continuous functions which possess a locally absolutely continuous $(r-1)$st derivative in $(-1,1)$, such that

$$\ddot{\mathbb{C}}^{(r)} = \{f \in \mathbb{C} : f^{(r-1)} \in AC_{loc}(-1,1) \text{ and } \|\phi^r f^{(r)}\| < \infty\}.$$

Now, we see appropriate to mention a concepts (co)convex and unconstrained polynomials.

**Definition 1.4 [25]** A subset $X$ of $\mathbb{R}^n$ is convex set if $[x,y] \subseteq X$, whenever $x, y \in X$. Equivalently, $X$ is convex if

$(1-\lambda)x + \lambda y \in X$, for all $x, y \in X$ and $\lambda \in [0,1]$.

The function $f$ is called convex of $X$ if

$$f((1-\lambda)x + \lambda y) \le (1-\lambda)f(x) + \lambda f(y), \text{ for all } x, y \in X \text{ and } \lambda \in [0,1].$$

**Definition 1.5 [21]** Let $\pi_n$ be the space of all algebraic polynomials of degree $\leq n-1$, and $\Delta^{(2)}$ be the set of all convex functions on $[-1,1]$. For $f \in \mathbb{C}([-1,1]) \cap \Delta^{(2)}$, the degree of best convex polynomial approximation of $f$ is

$$E_n^{(2)}(f) = \inf \{\|f - p_n\|, p_n \in \pi_n \cap \Delta^{(2)}\}.$$

**Definition 1.6 [8]** Let $Y_s = \{y_i\}_{i=1}^{s}$, $s \in \mathbb{N}$ be a partition of $[-1,1]$, that is, a collection of $s$ fixed points $y_i$ such that

$$y_{s+1} = -1 < y_s < \cdots < y_1 < 1 = y_o$$

and let $\Delta^{(2)}(Y_s)$ be the set of continuous functions on $[-1,1]$ that are convex downwards on the segment $[y_{i+1}, y_i]$ if $i$ is even and convex upwards on the same segment if $i$ is odd. The functions from $\Delta^{(2)}(Y_s)$ are called coconvex.

**Definition 1.7 [13]** Let $\Delta^{(2)}(Y_s)$ be the collection of all functions $f$ in $\mathbb{C}([-1,1])$ that change convexity at the points of the set $Y_s$, and are convex in $[y_s, 1]$. The degree of best CCP approximation of $f$ is defined by

$$E_n^{(2)}(f, Y_s) = \inf \{\|f - p_n\|, p_n \in \pi_n \cap \Delta^{(2)}(Y_s)\}.$$

**Definition 1.8 [3]** A domain $\mathbb{D}$ of convex polynomial $p_n$ of $\Delta^{(2)}$ is a subset of $X$ and $X \subseteq \mathbb{R}$, satisfying the following properties:

1) $\mathbb{D} \in \mathcal{K}^N$, where
$$\mathcal{K}^N = \{\mathbb{D}: \mathbb{D} \text{ is a compact subset of } X\}$$
is the class of all domain of convex polynomial,
2) there is the point $t \in X/\mathbb{D}$, such that
$$|p_n(t)| > \sup \{|p_n(x)|: x \in \mathbb{D}\}, \text{ and}$$
3) there is the function $f$ of $\Delta^{(2)}$, such that
$$\|f - p_n\| \leq \frac{c}{n^2} \omega_{2,2}^{\phi}\left(f'', \frac{1}{2}\right).$$

**Definition 1.9 [3]** A domain $\mathbb{D}$ of coconvex polynomial $p_n$ of $\Delta^{(2)}(Y_s)$ is a subset of $X$ and $X \subseteq \mathbb{R}$, satisfying the following properties:

1) $\mathbb{D} \in \mathcal{K}^N(Y_s)$, where
$$\mathcal{K}^N(Y_s) = \left\{ \begin{matrix} \mathbb{D}: \mathbb{D} \text{ is a compact subset of } X, \\ \text{and } p_n \text{ changes convexity at } \mathbb{D} \end{matrix} \right\}$$
is the class of all domain of coconvex polynomial,
2) $y_i$'s are inflection points, such that
$$|p_n(y_i)| \leq \frac{1}{2}, i = 1, \ldots, s, \text{ and}$$
3) there is the function $f$ of $\Delta^{(2)}(Y_s)$, such that
$$\|f - p_n\| \leq \frac{c}{n^2} \omega_{k,2}^{\phi}\left(f'', \frac{1}{n}\right).$$

**Definition 1.10 [18]** The polynomial $p_n$ is a constrained polynomial, if $p_n$ is shape preserving of function $f$, except that been UNP.

DeVore (1977) introduced estimates for the approximation of monotone functions by monotone polynomials (see [5]). Also, he given results of the order as for the unconstrained approximation by polynomials.

In 1992, DeVore et al. [4] proved for $f \in L_p$, $0 < p < 1$, and $k \in \mathbb{N}$, there exists an algebraic polynomial $p_n$ of degree $\leq n$ such that

$$\|f - p_n\|_p \leq C\omega_k^\phi\left(f, \frac{1}{n}\right)_p,$$

for $x \in I$, where $\omega_k^\phi$ is the usual DTMS.

In 1995, Kopotun [18] investigated the notion of degree of the best UNCP approximation presented in terms of $\omega_2^{\phi^\lambda}$. This approach of approximation does not depend on $\lambda$ and $0 \leq \lambda \leq 1$. In 1995, he proved uniform estimates for monotone and convex approximation of functions by algebraic polynomials in terms of the usual DTMS $\omega_k^\phi$.

**Definition 1.11 [26]** The partition $\tilde{T}_n = \{t_j\}_{j=0}^n$, where

$$t_j = t_{j,n} = \begin{cases} -\cos\left(\frac{j\Pi}{n}\right), & \text{if } 0 \leq j \leq n, \\ -1, & \text{if } j < 0, \end{cases}$$

and $t_j$'s as the knots of Chebyshev partition.

**Definition 1.12 [22]** Let $\Sigma_{k,n}$ be the collection of all continuous piecewise polynomials of degree $k - 1$, on the Chebyshev partition and let $\Sigma_{k,n}^1 \subseteq \Sigma_{k,n}$ be the subset of all continuously differentiable such functions. That is, if $\mathfrak{p} \in \Sigma_{k,n}$, then

$$\mathfrak{p}|_{I_j} = p_j, \quad j = 1, \ldots, n,$$

where $p_j \in \pi_{k-1}$, the collection of polynomials of degree $\leq k - 1$, and

$$p_j(x_j) = p_{j+1}(x_j), j = 1, \ldots, n - 1,$$

and if $\mathfrak{p} \in \Sigma_{k,n}^1$, then in addition,

$$p_j'(x_j) = p_{j+1}'(x_j), \quad j = 1, \ldots, n - 1.$$

The following theorems indicate several cases of different for inflection points.

**Theorem 1.13 [9] ($s \geq 2$)** Let $s \geq 2$. For every $\sigma > 0$, and each $Y_s \in \mathbb{Y}_s$ there exists $N(\sigma, Y_s)$ with the property that if $f \in \Delta^{(2)}(Y_s)$, then

$$\sup\left\{n^\sigma E_n^{(2)}(f, Y_s) : n > N(\sigma, Y_s)\right\} \leq c(\sigma, s) \sup\{n^\sigma E_n(f) : n \in \mathbb{N}\}.$$

**Theorem 1.14 [9] ($s = 1$)** For every $\sigma > 0$, $\sigma \neq 4$, $Y_1 \in \mathbb{Y}_1$ and $f \in \Delta^{(2)}(Y_1)$, we have

$$\sup\left\{n^\sigma E_n^{(2)}(f, Y_1) : n \in \mathbb{N}\right\} \leq c(\sigma) \sup\{n^\sigma E_n(f) : n \in \mathbb{N}\}.$$

**Theorem 1.15 [9] ($s = 1$)** There is an absolute constant $c$, such that for every $Y_1 \in \mathbb{Y}_1$ and a function $f \in \Delta^{(2)}(Y_1)$ the inequality

$$\sup\left\{n^4 E_n^{(2)}(f, Y_1): n > (1 - y_1^2)^{\frac{-1}{2}}\right\} \leq c \sup\{n^4 E_n(f): n \in \mathbb{N}\},$$

hold.

**Theorem 1.16 [9] ($s = 1$)** For every $Y_1 \in \mathbb{Y}_1$ there exists a function $f \in \Delta^{(2)}(Y_1)$ satisfying

$$\sup\{n^4 E_n(f): n \in \mathbb{N}\} = 1,$$

such that for each $m \in \mathbb{N}$, we have

$$n^4 E_m^{(2)}(f, Y_1) \geq c \ln\left(\frac{m}{1 + m^2 \phi(y_1)} - 1\right),$$

and

$$\sup\left\{n^4 E_n^{(2)}(f, Y_1): n \in \mathbb{N}\right\} \geq c \ln|\phi(y_1)|.$$

The paper structure is as follows: In Section 2, statements and definitions of the necessary and the main results are stated. Then, proof of our main theorems in Section 3. We show some applications of paper in Section 4.

## 2. Statements and Definitions of the Necessary and Main Results

In this section, we assume that the function $f$ is applying definition of the following.

**Definition 2.1 [2]** Let $\mathbb{D}$ be measurable set, $f: \mathbb{D} \to \mathbb{R}$ be a bounded function, and $\mathcal{L}_i: \mathbb{D} \to \mathbb{R}$ be nondecreasing function for $i \in \Lambda$. For a Lebesgue partition P of $\mathbb{D}$, put $\underline{LS}(f, P, \underline{\mathcal{L}}) = \sum_{j=1}^n \prod_{i \in \Lambda} m_j \, \mathcal{L}_i\left(\mu(\mathbb{D}_j)\right)$ and $\overline{LS}(f, P, \underline{\mathcal{L}}) = \sum_{j=1}^n \prod_{i \in \Lambda} M_j \, \mathcal{L}_i\left(\mu(\mathbb{D}_j)\right)$ where $\mu$ is a measure function of $\mathbb{D}$, $m_j = \inf\{f(x): x \in \mathbb{D}_j\}$, $M_j = \sup\{f(x): x \in \mathbb{D}_j\}$ and $\underline{\mathcal{L}} = \mathcal{L}_1, \mathcal{L}_2, \dots$. Also, $\mathcal{L}_i(x_j) - \mathcal{L}_i(x_{j-1}) > 0$, $\underline{LS}(f, P, \underline{\mathcal{L}}) \leq \overline{LS}(f, P, \underline{\mathcal{L}})$, $\prod_{i \in \Lambda} \underline{\int}_i^{\mathbb{D}} f \underline{d\mathcal{L}} = \sup\{\underline{LS}(f, \underline{\mathcal{L}})\}$ and $\prod_{i \in \Lambda} \overline{\int}_i^{\mathbb{D}} f \underline{d\mathcal{L}} = \inf\{\overline{LS}(f, \underline{\mathcal{L}})\}$ where $\underline{LS}(f, \underline{\mathcal{L}}) = \{\underline{LS}(f, P, \underline{\mathcal{L}}): P \text{ part of set } \mathbb{D}\}$ and $\overline{LS}(f, \underline{\mathcal{L}}) = \{\overline{LS}(f, P, \underline{\mathcal{L}}): P \text{ part of set } \mathbb{D}\}$. If $\prod_{i \in \Lambda} \underline{\int}_i^{\mathbb{D}} f \underline{d\mathcal{L}} = \prod_{i \in \Lambda} \overline{\int}_i^{\mathbb{D}} f \underline{d\mathcal{L}}$ where $\underline{d\mathcal{L}} = d\mathcal{L}_1 \times d\mathcal{L}_2 \times \dots$. Then $f$ is integral $\int_i$ according to $\mathcal{L}_i$ for $i \in \Lambda$.

**Lemma 2.2 [1]** If $f$ is a function of Lebesgue Stieltjes integral-i, then $vf$ is a function of Lebesgue Stieltjes integral-i, where $v > 0$ is real number, and

$$\prod_{i \in \Lambda} \int_i^{\mathbb{D}} vf \underline{d\mathcal{L}} = v \prod_{i \in \Lambda} \int_i^{\mathbb{D}} f \underline{d\mathcal{L}},$$

holds.

**Lemma 2.3 [1]** If the functions $f_1, f_2$ are integrable on the set $\mathbb{D}$ according to $\mathcal{L}_i$, for $i \in \Lambda$, then $f_1 + f_2$ is the function of integrable according to $\mathcal{L}_i$, for $i \in \Lambda$, such that

$$\prod_{i \in \Lambda} \int_i^{\mathbb{D}} (f_1 + f_2) d\underline{\mathcal{L}} = \prod_{i \in \Lambda} \int_i^{\mathbb{D}} f_1 d\underline{\mathcal{L}} + \prod_{i \in \Lambda} \int_i^{\mathbb{D}} f_2 d\underline{\mathcal{L}}.$$

**Remark 2.4 [1]** Let $\mathcal{I}_f$ be the class of all functions of integrable $f$ that satisfying Definition 2.1, i.e.,

$$\mathcal{I}_f = \{f : f \text{ is integrable function according to } \mathcal{L}_i, i \in \Lambda\}$$

$$= \left\{ f : \prod_{i \in \Lambda} \underline{\int_i}^{\mathbb{D}} f d\underline{\mathcal{L}} = \prod_{i \in \Lambda} \overline{\int_i}^{\mathbb{D}} f d\underline{\mathcal{L}} \right\}.$$

**Remark 2.5 [14]** Let $x_i \in \left[\frac{x_i + x^{\#}}{2}, \frac{x_i + x_*}{2}\right] \subseteq \theta_{\mathcal{N}}$, then we denote

$$x^{\#} = x_{j(i)+1}, \quad x_* = x_{j(i)-2}$$

where

$$\theta_{\mathcal{N}} = \theta_{\mathcal{N}}[-1,1] = \{x_i\}_{i=0}^{\mathcal{N}} = \{-1 = x_\circ \leq \cdots \leq x_{\mathcal{N}-1} \leq x_{\mathcal{N}} = 1\}$$

and

$$\|\theta_{\mathcal{N}}\| = \max_{0 \leq i \leq \mathcal{N}-1} \{x_{i+1} - x_i\}$$

the length of the largest interval in that partition.

**Definition 2.6 [1]** For $f \in \mathcal{I}_f$ and $r \in \mathbb{N}_o$, we define

$$\omega_{i,r}^{\phi}\left(f^{(r)}, \|\theta_{\mathcal{N}}\|, [-1,1]\right)_{w_{\alpha,\beta},p} = \sup\left\{ \left\|w_{\alpha,\beta} \phi^r \Delta_{h\phi}^i(f^{(r)}, x)\right\|_p, 0 < h \leq \|\theta_{\mathcal{N}}\| \right\},$$

where

$$\Delta_{h\phi}^i(f, x) = \prod_{i \in \Lambda} \int_i^{\mathbb{D}} f d\underline{\mathcal{L}}_{\phi}$$

and $\|\theta_{\mathcal{N}}\| < 2(i^{-1})$, $\mathcal{N} \geq 2$.

Now, we are ready to provide our next definition.

**Definition 2.7** For $\alpha, \beta \in J_p$ and $f \in \mathcal{I}_f$, we set

$$\mathbb{E}_n(f, w_{\alpha,\beta})_{\alpha,\beta,p} = \mathbb{E}_n(f)_{\alpha,\beta,p} = \inf\left\{ \|f - p_n\|_{\alpha,\beta,p}, \ p_n \in \pi_n \cap \mathcal{I}_f, \ f \in \Delta^{(2)}(Y_s) \cap \mathbb{L}_p^{\alpha,\beta} \cap \mathcal{I}_f \right\}$$

and

$$\mathcal{E}_n^{(2)}(f, w_{\alpha,\beta}, Y_s)_p = \inf\left\{ \|f - p_n\|_{\alpha,\beta,p}, \ p_n \in \pi_n \cap \Delta^{(2)}(Y_s) \cap \mathcal{I}_f, \ f \in \Delta^{(2)}(Y_s) \cap \mathbb{L}_p^{\alpha,\beta} \cap \mathcal{I}_f \right\}$$

respectively, denote the degree of best unconstrained and (co)convex polynomial approximation of $f$.

**Example 2.8** Consider the function

$$f(x) = \tan(\cos(\exp(x^4)))\qquad(1')$$

with a partition $Y_9 = \{y_i\}_{i=1}^{9}$ of $[-1,2]$, such that

$$y_{10} = -1 < 1.12 < \cdots < 1.92 < 2 = y_o.$$

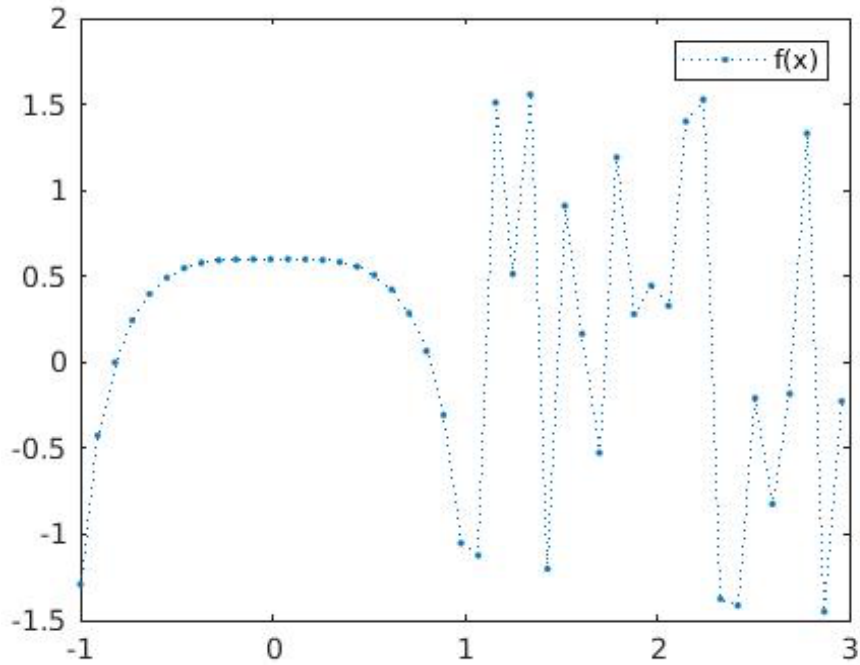

**Figure 1.** Coconvex function of $f \in \Delta^{(2)}(Y_9)$.

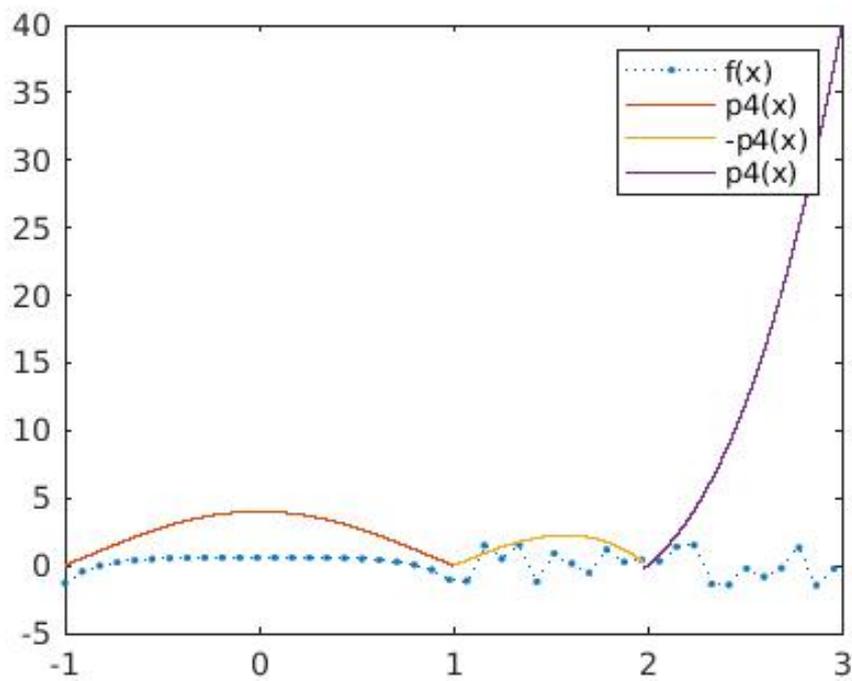

**Figure 2.** Show graph for CCP approximation.

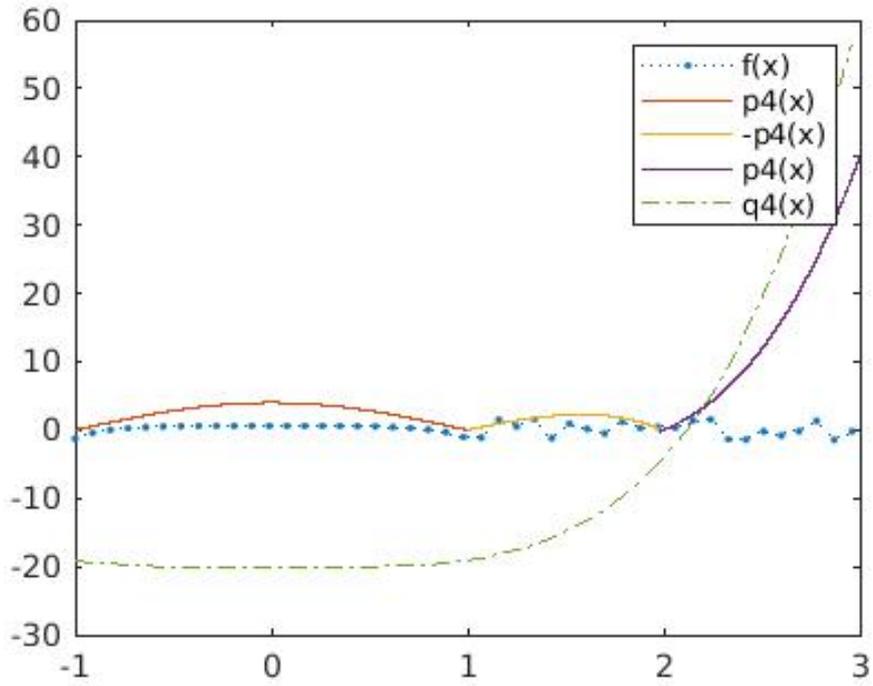

**Figure 3.** Show graph for COCUNP approximation.

These Figures 1-3 were verified with MATLAB. Suppose that $\alpha = \beta = 1$, then

$$\|f - p_n\|_{\alpha,\beta,p} = \|(1-x^2)(f-p_n)\|_p,$$

if $p = 1$, we have

$$= \int_{-1}^{2} |(1-x^2) \times \tan(\cos(\exp(x^4)))|dx - \int_{-1}^{2} |(1-x^2) \times p_4(x)|dx$$

$$= I(x) + \begin{cases} \int_{-1}^{2} (1-x^2)(x^4 - e^3)dx \\ \int_{-1}^{2} (1-x^2)(x+2)(x+1)(x-1)(x-2)dx. \end{cases}$$

Thus,

$$I(x) = I_o(x) + I_1(x)$$

$$I_o(x) = \left.\frac{\ln(\cos(\cos(\exp(x^4))))}{4x^3 \sin(\exp(x^4))}\right|_{-1}^{2} = -0.26$$

and

$$I_1(x) = \frac{\ln(\cos(\cos(\exp(x^4))))}{4x\,\sin(\exp(x^4))}\bigg|_{-1}^{2} = -0.34.$$

Finally, summary is drawn in the Table 1.

**Table 1.** Coconvex and unconstrained polynomial approximation.

| $y_i$ | Function | Polynomials | | | $\|f - p_n\|_{\alpha,\beta,p}$ | |
|---|---|---|---|---|---|---|
| | | $\pi_n$ | $\Delta^{(2)}(Y_s)$ | $\mathcal{J}_f$ | $\mathcal{E}_n^{(2)}$ | $\mathbb{E}_n$ |
| -1 | | | | | | |
| 1 | | | | | | |
| 1.12 | | $p_4(x) = x^4 - e^3$ | | + | | 11.218 |
| 1.28 | | | | | | |
| 1.36 | | | | | | |
| 1.44 | (1') | $p_4(x) =$ | | | | |
| 1.52 | | $\begin{cases}(x+2)(x+1)(x-1)(x-2); & \text{if } x \in [-1, 1.005], [1.981, 2] \\ -(x+2)(x+1)(x-1)(x-2); & \text{if } x \in [1.005, 1.981]\end{cases}$ | + | + | | −48.18 |
| 1.68 | | | | | | |
| 1.76 | | | | | | |
| 1.92 | | | | | | |
| 2 | | | | | | |

Now, it is appropriate to present the problem in this paper. Theorems 1.13-16 are invalid in $\mathbb{L}_p^{\alpha,\beta}$ space. However, these theorems were not ensuring derivation of degrees for best COCUNP approximation like the function (1'), such that $x \in [-1,2]/Q^c$ (see Figure 1).

**Theorem 2.9 [1]** For $r \in \mathbb{N}_o$, $\alpha, \beta \in J_p$, there is a constant $c = c(r, \alpha, \beta, p)$ such that if $f \in \Delta^{(2)} \cap \mathbb{L}_{p,r}^{\alpha,\beta}$, there, a number $\mathcal{N} = \mathcal{N}\left(f, \omega_{1,r}^\phi(f^{(r)}, \|\theta_\mathcal{N}\|, I)_{w_{\alpha,\beta},p}\right)$ for $n \geq \mathcal{N}$ and $\mathcal{S} \in \mathbb{S}(\tilde{T}_n, r+2) \cap \Delta^{(2)} \cap \mathbb{L}_{p,r}^{\alpha,\beta}$, such that

$$\left\|f^{(r)} - \mathcal{S}^{(r)}\right\|_{w_{\alpha,\beta},p} \leq c_{r,\alpha,\beta,p,\omega_{1,r}^\phi} \min\left\{\omega_{i,r}^\phi(f^{(r)}, \|\theta_\mathcal{N}\|, I_\alpha)_{w_{\alpha,\beta},p}, \omega_{i,r}^\phi(f^{(r)}, \|\theta_\mathcal{N}\|, I_\beta)_{w_{\alpha,\beta},p}\right\}$$

where

$$\left.\begin{aligned}\Delta_{h\phi,\alpha}^i(f^{(r)}, x) &= \int_1^\mathbb{D}\int_2^\mathbb{D}\cdots\int_i^\mathbb{D}\cdots f^{(r)}(x)\,d\mathcal{L}_{1t,\alpha}d\mathcal{L}_{2t,\alpha}\cdots d\mathcal{L}_{it,\alpha}\cdots = \prod_{i\in\Lambda}\int_i^\mathbb{D} f^{(r)}\,\underline{d\mathcal{L}_{t\phi,\alpha}}\\ \Delta_{h\phi,\beta}^i(f^{(r)}, x) &= \int_1^\mathbb{D}\int_2^\mathbb{D}\cdots\int_i^\mathbb{D}\cdots f^{(r)}(x)\,d\mathcal{L}_{1t,\beta}d\mathcal{L}_{2t,\beta}\cdots d\mathcal{L}_{it,\beta}\cdots = \prod_{i\in\Lambda}\int_i^\mathbb{D} f^{(r)}\,\underline{d\mathcal{L}_{t\phi,\beta}}\end{aligned}\right\} \quad (1)$$

Moreover, if $r, \alpha, \beta = 0$, then

$$\|f - \mathcal{S}\|_p \leq c\left(\omega_1^\phi\right)\omega_i^\phi(f, \|\theta_\mathcal{N}\|, I)_p.$$

In particular,

$$\left\|f^{(r)} - \mathcal{S}^{(r)}\right\|_{w_{\alpha,\beta},p} \leq c_r\,\omega_{1,r}^\phi(f^{(r)}, \|\theta_\mathcal{N}\|, I)_{w_{\alpha,\beta},p}.$$

**Theorem 2.10 [1]** Let $\Delta^k$ be the space of all $k$-monotone functions. If $f \in \Delta^k \cap \mathbb{L}_{p,r}^{\alpha,\beta}$ is such that $f^{(r)}(x) = p_n^{(r)}(x)$, where $p_n \in \pi_n \cap \Delta^k$, $N \geq k \geq 2$ and $s \in \mathbb{S}(\tilde{T}_n, r+2) \cap \Delta^k \cap \mathbb{L}_{p,r}^{\alpha,\beta}$. Then

$$\|f - s\|_{w_{\alpha,\beta},p} \leq c(f, p, k, \alpha, \beta, x_\star, x^\#)\omega_{i,r}^\phi(f, \|\theta_\mathcal{N}\|, I)_{w_{\alpha,\beta},p}.$$

In particular, if $f$ is a convex function and $p_n$ is a convex polynomial or piecewise convex polynomial, then

$$\|f - s\|_{w_{\alpha,\beta},p} \leq c_k \omega_{i,r}^{\phi}(f, \|\theta_{\mathcal{N}}\|, I)_{w_{\alpha,\beta},p}.$$

Now, the main results are ready for take positions in this paper.

**Theorem 2.11** Let $\sigma, m, n \in \mathbb{N}$, $\sigma \neq 4$, $s \in \mathbb{N}_o$ and $\alpha, \beta \in J_p$. If $f \in \Delta^{(2)}(Y_s) \cap \mathbb{L}_p^{\alpha,\beta} \cap \mathcal{I}_f$, then

$$\sup\left\{n^\sigma \mathcal{E}_n^{(2)}(f, w_{\alpha,\beta}, Y_s)_p : n \geq m\right\} \leq c \sup\{n^\sigma \mathbb{E}_n(f)_{\alpha,\beta,p} : n \in \mathbb{N}\}. \tag{2}$$

In particular, suppose that $Y_s \in \mathbb{Y}_s$ and $s \geq 1$. Then

$$\mathcal{E}_n^{(2)}(f, w_{\alpha,\beta}, Y_s)_p \leq cn^{-\sigma} \omega_{i,r}^{\phi}(f^{(r)}, \|\theta_{\mathcal{N}}\|, I)_{\alpha,\beta,p}, n \geq \|\theta_{\mathcal{N}}\|.$$

**Theorem 2.12** If $\sigma, n \in \mathbb{N}$, $\sigma = 4$, $I = [-1,1]$, $Y_1 \in \mathbb{Y}_1$, $\alpha, \beta \in J_p$ and $c = c(\alpha, \beta, p)$. If $f \in \Delta^{(2)}(Y_s) \cap \mathbb{L}_p^{\alpha,\beta} \cap \mathcal{I}_f$, then,

$$n^\sigma \mathbb{E}_n(f)_{\alpha,\beta,p} \leq c,$$

and

$$n^{-\eta} \|f\|_{\alpha,\beta,p} \leq c \times \mathcal{E}_n^{(2)}(f, w_{\alpha,\beta}, Y_1)_p, \tag{3}$$

where $\eta = i + \text{dist}(I, \{-1,1\})$.

### 3. Proof of Main Theorems

Before proving our results, we need the following remark.

**Remark 3.1** If $f$ in $\mathcal{I}_f$ is a function of Lebesgue Stieltjes integral-i, and $f$ is differentiable function, therefore,

$$f' = \frac{df}{dx} = \frac{d}{dx}\left(\int_0^x \frac{df(u)}{d\mathcal{G}_{1,\mu,\mathbb{D}_o}} d\mathcal{G}_{1,\mu,\mathbb{D}_o}\right)$$

$$= \frac{d}{dx}\left(\int_0^x \int_0^x \frac{d^2 f(u)}{d\mathcal{G}_{1,\mu,\mathbb{D}_o} \times d\mathcal{G}_{2,\mu,\mathbb{D}_o}} d\mathcal{G}_{1,\mu,\mathbb{D}_o} \times d\mathcal{G}_{2,\mu,\mathbb{D}_o}\right)$$

$$= \frac{d}{dx}\left(\int_0^x \int_0^x \cdots \int_0^x \cdots \frac{d^i f(u)}{d\mathcal{G}_{1,\mu,\mathbb{D}_o} \times d\mathcal{G}_{2,\mu,\mathbb{D}_o} \times \cdots \times d\mathcal{G}_{i,\mu,\mathbb{D}_o} \times \cdots} d\mathcal{G}_{1,\mu,\mathbb{D}_o} \times d\mathcal{G}_{2,\mu,\mathbb{D}_o} \times \cdots \times d\mathcal{G}_{i,\mu,\mathbb{D}_o} \times \cdots\right)$$

$$= \frac{d}{dx}\left(\prod_{i \in \Lambda} \int_i^{I_x} f^{(i)}(u) \, \underline{d\mathcal{G}_{\mu,\mathbb{D}_o}}\right), x \in I_x = [0, x] \subseteq \mathbb{D}_o, u \in \mathbb{D}_o \text{ and } \mathcal{G}_{\mu,\mathbb{D}_o} = \mathcal{G}(\mu(\mathbb{D}_o))$$

$$= \frac{d}{dx}\left(\prod_{i\in\Lambda}\int_i^{I_x} f^{(i)}\,\underline{d\mathcal{G}_{\mu,\mathbb{D}_o}}\right)$$

$$= \prod_{i\in\Lambda}\int_i^{I_x} f_x^{(i+1)}\,\underline{d\mathcal{G}_{\mu,\mathbb{D}_o}}\,.$$

$$f_x^{(i+1)} = \frac{d^i}{d\mathcal{G}_{i,\mu,\mathbb{D}_o}^i}f' = \frac{d}{dx}\left(\frac{d^i f}{d\mathcal{G}_{i,\mu,\mathbb{D}_o}^i}\right) = \frac{d^{i+1}f_x}{dx\times d\mathcal{G}_{i,\mu,\mathbb{D}_o}^i}\,.$$

**3.2 Proof of Theorem 2.11** Let $f$ be a fixed shape preserving (convexity) in ($I_1 = [y_2, 1]$ and $I_n = [-1, y_{n-2}]$). Now, if $J = [a,b] \subseteq [-1,1]$, we should use $/J/ = \frac{|J|}{\phi\left(\frac{\|\theta_N\|}{2}\right)}$ and $\phi(x) = \sqrt{1-x^2}$ (see [23]).

By using Definition 2.1, technique of Lemmas 2.2, 2.3 and inequality above, (with ease, it becomes in our inventory the following conclusion).

$$\left|I_j\right|^{\ell}\omega_m\!\left(f^{(r)},|I_j|,I_j\right)_{\alpha,\beta,p} = \left|I_j\right|^{\ell}\times\left\{\sup\left\|w_{\alpha,\beta}\phi^{\ell}\Delta_h^m(f^{(r)},.)\right\|_p, 0 < h \le \frac{|I_j|}{i}\phi\right\}$$

where $m = i + r - \ell$ and $0 \le \ell \le r \le i$,

$$\left|I_j\right|^{\ell}\omega_m\!\left(f^{(r)},|I_j|,I_j\right)_{\alpha,\beta,p} = \left|I_j\right|^{\ell}\times$$

$$\left\{\left|\sup\left(w_{\alpha,\beta}\times\sup\left(\sum_{j=1}^n\prod_{m\in\Lambda}\inf_j\!\left(f^{(r)}(x)\right)\mathcal{L}_{mt}\!\left(\mu(\mathbb{D}_j)\right)\right)\times\mu(\mathbb{D}_j)\right)\right|, 0 < h \le \frac{|I_j|}{i}\phi\right\},$$

since $m = i + r - \ell$, then $i - \ell < i$ and $i + r - \ell \le 2i$, whence, $\frac{2i}{m} \ge 1$,

$$\left|I_j\right|^{\ell}\omega_m\!\left(f^{(r)},|I_j|,I_j\right)_{\alpha,\beta,p} \le \left\{\left(\left|\sup L(w_{\alpha,\beta}\times\sup LS(f^{(r)},\mathcal{L}_t))\right|^p\right)^{\frac{1}{p}}, 0 < h \le \frac{|I_j|}{i}\right\}$$

$$\le cm^{-1}\left\{\left(\left|\sup\left(\sum_{j=1}^n\prod_{i\in\Lambda}\inf_j\!\left(f^{(r)}(x)\right)\mathcal{L}_{it}\!\left(\mu(\mathbb{D}_j)\right)\right)\right|^p\right)^{\frac{1}{p}}, 0 < h \le \|\theta_N\|\right\}$$

$$\le cm^{-1}\omega_{i,r}^{\phi}(f^{(r)},\|\theta_N\|)_{\alpha,\beta,p}.$$

Assume that $\mathfrak{p}_n \in \Sigma_{i+5,n}(Y_s)\cap\Delta^{(2)}(Y_s)$ and the estimate

$$\omega_{i+5}^{\phi}(\mathfrak{p}_n, n^{-1})_{\alpha,\beta,p} = \left\{\sup\left\|w_{\alpha,\beta}\phi^r\Delta_{h\phi}^{i+5}(\mathfrak{p}_n,.)\right\|_p, 0 < h \le \|\theta_N\|\right\}$$

$$\le c\left\{\sup\left\|w_{\alpha,\beta}\phi^r\Delta_{h\phi}^{i+5}(\mathfrak{p}_n - f + f,.)\right\|_p, 0 < h \le \|\theta_N\|\right\}. \tag{4}$$

Now,

$$\Delta_{h\phi}^{i+5}(\mathfrak{p}_n,.) = \Delta_{h\phi}^{i+5}(\mathfrak{p}_n - f + f,.) = \prod_{i\in\Lambda}\int_{i+5}^{\mathbb{D}}(\mathfrak{p}_n - f + f)\,\underline{d\mathcal{L}_{\phi}}$$

$$= \prod_{i \in \Lambda} \int_{i+5}^{\mathbb{D}} (\mathfrak{p}_n - f) \, d\mathcal{L}_\phi + \prod_{i \in \Lambda} \int_{i+5}^{\mathbb{D}} f \, d\mathcal{L}_\phi = \Delta_{h\phi}^{i+5}(\mathfrak{p}_n - f,.) + \Delta_{h\phi}^{i+5}(f,.). \tag{5}$$

Then, the proofs (4) and (5) an immediately give the following:

$$\omega_{i+5}^{\phi}(\mathfrak{p}_n, n^{-1})_{\alpha,\beta,p} \leq c\|f - \mathfrak{p}_n\|_{\alpha,\beta,p} + c\left\{\sup\|w_{\alpha,\beta}\phi^r \Delta_{h\phi}^{i+5}(f^{(5)},.)\|_p, 0 < h \leq \|\theta_\mathcal{N}\|\right\}$$

$$\leq c\|f - \mathfrak{p}_n\|_{\alpha,\beta,p} + cm^{-1}\omega_{i+5}^{\phi}(f^{(5)}, n^{-1})_{\alpha,\beta,p}.$$

If $(\ell, r) = (3,5)$ and $(\ell, r) = (2,5)$ are implies, respectively,

$$|I_j|^3 \omega_{i+2}(f^{(3)}, |I_j|, I_j)_{\alpha,\beta,p} \leq cm^{-1}\omega_{i,5}^{\phi}(f^{(5)}, \|\theta_\mathcal{N}\|)_{\alpha,\beta,p}, 1 < j < n,$$

and

$$|I_j|^2 \omega_{i+3}(f'', |I_j|, I_j)_{\alpha,\beta,p} \leq cm^{-1}\omega_{i,5}^{\phi}(f^{(5)}, \|\theta_\mathcal{N}\|)_{\alpha,\beta,p}, 1 \leq j \leq n.$$

Let $O_i = [a_i, b_i]$ be the sets, $1 \leq i \leq s$, are all disjoint and do not contain the endpoints of $[-1,1]$. We have chosen an arbitrary constant $N(Y_s)$ depend on $O_i$, for $n \geq N$.

Next, if $O_i \cap I_i = \emptyset$, then $f$ does not change its convexity on $I_i$. Let $k \geq 1$, and $a_i < c < a_i + h$ be fixed-large enough. The inequality

$$f''(x)(x - c) \geq 0, \quad a_i \leq x \leq a_i + h,$$

is an immediate result of $f \in \Delta^{(2)}(Y_s) \cap \mathbb{L}_p^{\alpha,\beta}$.

**PART I. (Co)convex Polynomial.** If $p_n \in \pi_n \cap \Delta^{(2)}(Y_s)$, $s \geq 0$ such that

$$p_n(x)(x - c) \geq 0, \quad a \leq x \leq a + h,$$

then, there exists the sequence $f_n$ and

$$p_n(x) = \begin{cases} \dfrac{f'(a)}{f(a)} & ; \quad \text{if } f_n(a+h) + f'(a) \leq f'(a+h), \\ \dfrac{f'_n(a+h) + f'(a+h)}{f_n(a+h) + f(a+h)} & ; \quad \text{otherwise.} \end{cases}$$

Now, we are willing to derive the following estimates based on the above outcomes.

$$\|f - p_n\|_{\alpha,\beta,p} = \|w_{\alpha,\beta}(f - p_n)\|_p \leq \|w_{\alpha,\beta}(f - p_n)\|_{L_p(\mathbb{D}_o)}, \mathbb{D}_o = [a, a+h].$$

**Case I.** If

$$f_n(a+h) + f'(a) \leq f'(a+h), \tag{6}$$

then,

$$\|f - p_n\|_{\alpha,\beta,p} \leq \left\|w_{\alpha,\beta}\left(f - \dfrac{f'(a)}{f(a)}\right)\right\|_{L_p(\mathbb{D}_o)}$$

$$\leq \left( \int_{\mathbb{D}_o} \left| w_{\alpha,\beta} \left( f - \frac{f'(a)}{f(a)} \right) \right|^p dx \right)^{\frac{1}{p}}, \quad x \in \mathbb{D}_o$$

$$\leq c_1 \left( \int_{\mathbb{D}_o} \left| w_{\alpha,\beta} \left( f - \int_{\mathbb{D}_o} \frac{f'(a)}{f(a)} d\mathcal{G}_{\mu,\mathbb{D}_o} \right) \right|^p dx \right)^{\frac{1}{p}}, \text{ where } \mathcal{G}_{\mu,\mathbb{D}_o} = \mathcal{G}(\mu(\mathbb{D}_o)) \text{ and } c_1 = c(\mathcal{G},\mu)^{-1}$$

$$\leq c_1 \left[ \left( \int_{\mathbb{D}_o} |w_{\alpha,\beta} f|^p dx \right)^{\frac{1}{p}} + \left( \int_{\mathbb{D}_o} \left| w_{\alpha,\beta} \times \int_{\mathbb{D}_o} \frac{f'(a)}{f(a)} d\mathcal{G}_{\mu,\mathbb{D}_o} \right|^p dx \right)^{\frac{1}{p}} \right]$$

$$\leq c_1 \left[ \left( |\sup L(w_{\alpha,\beta} f)|^p \right)^{\frac{1}{p}} + \left( \left| \sup L \left( w_{\alpha,\beta} \times \sup LS \left( \frac{f'(a)}{f(a)}, \mathcal{G}_{\mu,\mathbb{D}_o} \right) \right) \right|^p \right)^{\frac{1}{p}} \right]$$

$$\leq c_1 \left[ \left( |\sup L(w_{\alpha,\beta} f)|^p \right)^{\frac{1}{p}} + \left( \left| \sup LS \left( \sup L \left( w_{\alpha,\beta} \times \frac{f'(a)}{f(a)} \right), \mathcal{G}_{\mu,\mathbb{D}_o} \right) \right|^p \right)^{\frac{1}{p}} \right], \quad (7)$$

the following is an immediate inequality of $\int_{\mathbb{D}_o} \frac{f'(a)}{f(a)} dx = \ln f(a), x \in \mathbb{D}_o$.

$$\|f - p_n\|_{\alpha,\beta,p} \leq c(\mathcal{G},\mu,a)^{-1} \left[ \left( |\sup L(w_{\alpha,\beta} f)|^p \right)^{\frac{1}{p}} + \left( |\sup LS(w_{\alpha,\beta} \times \ln f(a), \mathcal{G}_{\mu,\mathbb{D}_o})|^p \right)^{\frac{1}{p}} \right],$$

from (1), then

$$\|f - p_n\|_{\alpha,\beta,p} \leq c(\mathcal{G},\mu,a)^{-1} \left[ \left\{ \sup \|w_{\alpha,\beta} \Delta_h^1(f,.)\|_p, \ 0 < h \leq t \right\} + \right.$$

$$\left\{ \sup \|w_{\alpha,\beta} \Delta_h^1(f,.)\|_p, 0 < h \leq t \right\} \right]$$

$$\leq 2 \times c(\mathcal{G},\mu,a)^{-1} \omega_1(f,t)_{\alpha,\beta,p}$$

$$\leq 2 \times c m^{-1} \omega_1^\phi(f,t)_{\alpha,\beta,p}.$$

**Case I. A. ($s = 0$)** We write $q_{n,r}(f,. |x_o, \ldots, x_{k-1}, x_*, x^\#)$ is defined in proof of Theorem 2.10 (see [1]), and $s(x)$ replace to $\mathfrak{p}(x)$. For any $v \in I_1$, and $I_n$ are preserving its convexity with the function $f$. Virtue of [1, eq. (10)], $r = 0$, we get

$$\|f - \mathfrak{p}\|_{\alpha,\beta,p} \leq c_{k,v} \omega_i(f, |I_1|, I_1)_{\alpha,\beta,p} \leq c_{k,v} \omega_i^\phi(f, \|\theta_N\|)_{\alpha,\beta,p}.$$

**Case I. B. ($s \geq 1$)** From (7), and Remark 3.1, we have

$$\|f - p_n\|_{\alpha,\beta,p} \leq \left\| w_{\alpha,\beta} \left( f - \frac{f'(a)}{f(a)} \right) \right\|_{L_p(\mathbb{D}_o)}$$

$$\leq c_1 \left[ \left( |\sup L(w_{\alpha,\beta} f)|^p \right)^{\frac{1}{p}} + \left( \left| \sup LS \left( \sup L \left( w_{\alpha,\beta} \times \left[ \frac{f'(a)}{f(a)} = F \right] \right), \mathcal{G}_{\mu,\mathbb{D}_o} \right) \right|^p \right)^{\frac{1}{p}} \right]$$

$$\leq c_1\left[\left(|\sup L(w_{\alpha,\beta}f)|^p\right)^{\frac{1}{p}} + \left(|\sup LS(\sup L(w_{\alpha,\beta} \times \inf L(F')), \mathcal{G}_{\mu,\mathbb{D}_o})|^p\right)^{\frac{1}{p}}\right],$$

since

$$F' = \frac{dF}{dx} = \prod_{i\in\Lambda}\int_i^{I_x} F_x^{(i+1)}\, d\mathcal{G}_{\mu,\mathbb{D}_o} = \sup LS\left(F_x^{(i+1)}, \mathcal{G}_{\mu,\mathbb{D}_o}\right),$$

then,

$$\|f - p_n\|_{\alpha,\beta,p} \leq c_1\left[\left(|\sup L(w_{\alpha,\beta}f)|^p\right)^{\frac{1}{p}} + \right.$$

$$\left.\left(\left|\sup LS\left(\sup L\left(w_{\alpha,\beta} \times \inf L\left(\sup LS\left(F_x^{(i+1)}, \mathcal{G}_{\mu,\mathbb{D}_o}\right)\right)\right), \mathcal{G}_{\mu,\mathbb{D}_o}\right)\right|^p\right)^{\frac{1}{p}}\right]$$

$$\leq c_1\left[\left(|\sup L(w_{\alpha,\beta}f)|^p\right)^{\frac{1}{p}} + \right.$$

$$\left.\left(\left|\sup LS\left(\left(\sup L\left(w_{\alpha,\beta} \times \sup LS\left(F^{(i)}, \mathcal{G}_{\mu,\mathbb{D}_o}\right)\right)\right), \mathcal{G}_{\mu,\mathbb{D}_o}\right)\right|^p\right)^{\frac{1}{p}}\right]$$

$$\leq c_1\left[\left(|\sup L(w_{\alpha,\beta}f)|^p\right)^{\frac{1}{p}} + \left(\left|\sup LS\left(\sup L\left(w_{\alpha,\beta} \times \left(\frac{f'(a)}{f(a)}\right)^{(i-1)}\right), \mathcal{G}_{\mu,\mathbb{D}_o}\right)\right|^p\right)^{\frac{1}{p}}\right],$$

by (1), we have

$$\|f - p_n\|_{\alpha,\beta,p} \leq 2c_1\left[\left(|\sup L(w_{\alpha,\beta}f)|^p\right)^{\frac{1}{p}} + \sup\left\{\|w_{\alpha,\beta}\phi^{i-1}\Delta_{h\phi}^{i-1}(f^{(i-1)},.)\|_p, 0 < h \leq \|\theta_\mathcal{N}\|\right\}\right]$$

$$\leq 2c_1\left[\|f\|_{\alpha,\beta,p} + \omega_{i-1,i-1}^{\phi}(f^{(i-1)}, t)_{\alpha,\beta,p}\right].$$

**Case II.** If (6) is invalid, hence, create estimates by the second part of the polynomial $p_n$, i.e.,

$$\|f - p_n\|_{\alpha,\beta,p} \leq \left\|w_{\alpha,\beta}\left(f - \left(\frac{f_n'(a+h) + f'(a+h)}{f_n(a+h) + f(a+h)}\right)\right)\right\|_{L_p(\mathbb{D}_o)}$$

$$\leq c(a+h, n, f)\omega_1^{\phi}(f, t)_{\alpha,\beta,p}.$$

**PART II.** *Unconstrained Polynomial.* Assume that, $(-1)^{i-s}\bigl(f(x) - q_{i-1}(f,.\,|y_1,\ldots,y_s)\bigr) \geq 0$, $x \in (y_i, y_{i+1})$, but, without absence of generality.

$$\mathbb{E}_n(f)_{\alpha,\beta,p} \leq \{\|f - \hat{p}_n\|_{\alpha,\beta,p}, \hat{p}_n \in \pi_n\}$$

$$\leq \|w_{\alpha,\beta}(f - \hat{p}_n)\|_p \leq c(f, p, i, \alpha, \beta, x_\star, x^\#)\omega_{i,r}^{\phi}(f, \|\theta_\mathcal{N}\|, I)_{\alpha,\beta,p}.$$

Hence, the theorem is proved. □

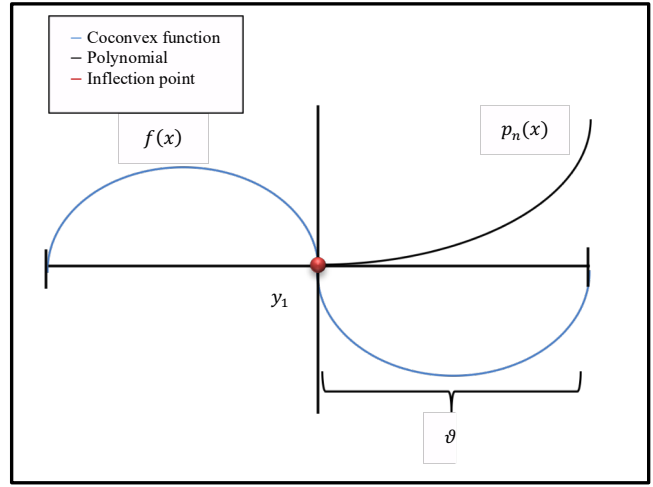

**Figure 4.** The degree of UNP approximation of coconvex function $f$ and $y_1 \in Y_1$.

**3.3 Proof of Theorem 2.12** Assuming $f \in \Delta^{(2)}(Y_1) \cap \mathbb{L}_p^{\alpha,\beta} \cap \mathcal{I}_f$, such that

$$f''(x)(x - y_1) \geq 0, \quad y_1 \in Y_1.$$

Let

$$(x - \vartheta)_+^{n-1} = \max\{0, (x - \vartheta)^{n-1}\}$$

and

$$p_n(x) = \frac{(1 - \vartheta)^{1-n}}{n}(x - \vartheta)_+^{n-1}, \quad n \in \mathbb{N}, \vartheta \in (0,1) \text{ and } x \in I,$$

be an arbitrary algebraic polynomial of $\leq n - 1$. Next, by Definitions 2.1, 2.7 and technique Lammas 2.2, 2.3, then

$$\mathbb{E}_n(f)_{\alpha,\beta,p} = \inf\left\{\|f - p_n\|_{\alpha,\beta,p}, \ p_n \in \pi_n \cap \mathcal{I}_f, \ f \in \Delta^{(2)}(Y_1) \cap \mathbb{L}_p^{\alpha,\beta} \cap \mathcal{I}_f\right\}$$

$$\leq \|f - p_n\|_{\alpha,\beta,p} = \|w_{\alpha,\beta}(f - p_n)\|_p$$

$$= \left(\int_{-1}^{1}|w_{\alpha,\beta}(f - p_n)|^p \, dx\right)^{\frac{1}{p}} = \left(\int_{-1}^{1}|w_{\alpha,\beta}f - w_{\alpha,\beta}p_n|^p \, dx\right)^{\frac{1}{p}}$$

$$= \left(\int_{-1}^{1}\left|w_{\alpha,\beta}f - w_{\alpha,\beta}\left[\frac{(1-\vartheta)^{1-n}}{n}(x-\vartheta)_+^{n-1}\right]\right|^p dx\right)^{\frac{1}{p}}$$

$$\leq \left(\int_{-1}^{1}|w_{\alpha,\beta}f|^p \, dx\right)^{\frac{1}{p}} + \left[\frac{(1-\vartheta)^{1-n}}{n}\right]\left(\int_{-1}^{1}|w_{\alpha,\beta}(x-\vartheta)_+^{n-1}|^p \, dx\right)^{\frac{1}{p}}$$

$$\leq \left(\int_{-1}^{1}|w_{\alpha,\beta}f|^{p}\,dx\right)^{\frac{1}{p}}+\left[\frac{(1-\vartheta)^{1-n}}{n}\right]\times$$

$$\left[\left(\int_{-1}^{0}|w_{\alpha,\beta}(x-\vartheta)_{+}^{n-1}|^{p}\,dx\right)^{\frac{1}{p}}+\left(\int_{0}^{1}|w_{\alpha,\beta}(x-\vartheta)_{+}^{n-1}|^{p}\,dx\right)^{\frac{1}{p}}\right].$$

From Figure 4, then

$$\|w_{\alpha,\beta}(f-p_{n})\|_{p}^{p}\leq\int_{-1}^{1}|w_{\alpha,\beta}f|^{p}\,dx+\left[\frac{(1-\vartheta)^{1-n}}{n}\right]\times$$

$$\left[\int_{-1}^{0}w_{\alpha,\beta}|(x-\vartheta)_{+}^{n-1}|^{p}\,dx+\int_{0}^{1}w_{\alpha,\beta}|(x-\vartheta)_{+}^{n-1}|^{p}\,dx\right]$$

$$\leq\int_{-1}^{1}|w_{\alpha,\beta}f|^{p}\,dx+c\times\left[\frac{(1-\vartheta)^{1-n}}{n}\right]\times\left[\int_{0}^{1}w_{\alpha,\beta}|(x-\vartheta)_{+}^{n-1}|^{p}\,dx\right], \quad \text{if } n=5$$

$$\leq\int_{-1}^{1}|w_{\alpha,\beta}f|^{p}\,dx+c\times\left[\frac{(1-\vartheta)^{-4}}{5}\right]\times\begin{cases}\left[\int_{0}^{1}w_{\alpha,\beta}|(x-\vartheta)_{+}^{4}|^{p}\,dx\right] & ;\text{if } x>\vartheta,\\ 0; & \text{if } x\leq\vartheta.\end{cases}$$

Therefore,

$$m^{4}\times\mathbb{E}_{n}(f)_{\alpha,\beta,p}\leq\left(\int_{-1}^{1}|w_{\alpha,\beta}f|^{p}\,dx\right)^{\frac{1}{p}}+c\times\begin{cases}\left[\left(\int_{0}^{1}w_{\alpha,\beta}|(x-\vartheta)_{+}^{4}|^{p}\,dx\right)^{\frac{1}{p}}\right] & ;\text{if } x>\vartheta,\\ 0; & \text{if } x\leq\vartheta,\end{cases}$$

where $m\leq n$.

Let's check the reliability of (3).

$$\|f\|_{\alpha,\beta,p}=\left(\int_{-1}^{1}|w_{\alpha,\beta}f|^{p}\,dx\right)^{\frac{1}{p}}=\left(\int_{-1}^{1}|w_{\alpha,\beta}(f-p_{n}+p_{n})|^{p}\,dx\right)^{\frac{1}{p}},$$

such that $p_{n}\in\Delta^{(2)}(Y_{1})\cap\pi_{n}\cap\mathcal{I}_{f}$

$$\|f\|_{\alpha,\beta,p}\leq\left(\int_{-1}^{1}|w_{\alpha,\beta}(f-p_{n})|^{p}\,dx\right)^{\frac{1}{p}}+\left(\int_{-1}^{1}|w_{\alpha,\beta}\times p_{n}|^{p}\,dx\right)^{\frac{1}{p}}$$

$$\leq c\left[\mathcal{E}_{n}^{(2)}(f,w_{\alpha,\beta},Y_{1})_{p}+\|p_{n}\|_{\alpha,\beta,p}\right], \quad n\in\mathbb{N}$$

$$\leq cn^{\eta}\mathcal{E}_{n}^{(2)}(f,w_{\alpha,\beta},Y_{1})_{p}.$$

Therefore,

$$n^{-\eta}\|f\|_{\alpha,\beta,p} \leq c \times \mathcal{E}_n^{(2)}(f, w_{\alpha,\beta}, Y_1)_p. \qquad \square$$

## 4. Some Applications

In this section, we will present some applications which type overlapping with COCUNP approximation like a Korovkin type approximation theorem.

Fejer operators are

$$F_n(f;x) = \frac{a_o}{2} + \sum_{k=1}^{n} \frac{n-k}{n}\left(a_k \frac{kx - x^\#}{x_\star - x_i} + b_k \frac{kx - x_\star}{x_i - x^\#}\right).$$

In 2014, Al-Muhja [2] defined $T_n$ by

$$T_n(f;x) = \frac{a_o}{2} + \sum_{k=1}^{n} \lambda_k^{(n)}\left(a_k \frac{kx - x^\#}{x_\star - x_i} + b_k \frac{kx - x_\star}{x_i - x^\#}\right),$$

where $\left\{\lambda_k^{(n)}\right\}_{k=1}^{n}$ is a matrix of real numbers, $n = 1,2,\ldots$ and also $a_k$ and $b_k$ are Fourier coefficients.

**Theorem 4.1 [2]** If a sequence $G_n(f)$ is positive linear functional and bounded on $\mathbb{C}(S)$, $f$ is bounded measurable function to $S$. Then, there exists nondecreasing function to $S$, such that

$$st_A - \lim_{\mu(S) \to 0} \left(\sup_n G_n(f) - f\right) = 0.$$

Assume that $A = (a_k)$ is a nonnegative regular summability matrix and

$$st_A - \lim_{t \to 0} \sum_{v=1}^{3} \omega_{i-1,i-1}^{\phi}\left(f_v^{(i-1)}, t\right)_{\alpha,\beta,p} = 0. \qquad (8)$$

**Theorem 4.2** Suppose that $S \subseteq \mathbb{R}$ is a Lebesgue measurable and $\mu(S) < \infty$. If a sequence $f_n(x)$ is Lebesgue measurable function and $f_n(x)$ is finite a. e. as $n$ is finite for $x \in S$. Then,

$$st_A - \lim_{n \to \infty} \mathcal{E}_n^{(2)}(f_n, w_{\alpha,\beta}, Y_s)_p \sim st_A - \lim_{n \to \infty} \mathbb{E}_n(f_n)_{\alpha,\beta,p}, \qquad (9)$$

where $f_n \in \Delta^{(2)}(Y_s) \cap \mathbb{L}_p^{\alpha,\beta} \cap \mathcal{I}_f$.

Proof. We let (2) be hold. Assuming

$$\mathcal{E}_n^{(2)}(f_n, w_{\alpha,\beta}, Y_s)_p \leq \mathbb{E}_n(f_n)_{\alpha,\beta,p}, \qquad (10)$$

then, (10) implies $f_n \in \Delta^{(2)}(Y_s) \cap \mathbb{L}_p^{\alpha,\beta} \cap \mathcal{I}_f$ (see Theorem 2.11), and put $S$ is Lebesgue measurable. Now,

$$\mathcal{E}_n^{(2)}(f_n, w_{\alpha,\beta}, Y_s)_p = \inf\left\{\|f_n - p_n\|_{\alpha,\beta,p}, \ p_n \in \Delta^{(2)}(Y_s) \cap \pi_n \cap \mathcal{I}_f, \ f_n \in \Delta^{(2)}(Y_s) \cap \mathbb{L}_p^{\alpha,\beta} \cap \mathcal{I}_f\right\}$$

$$\leq \|f_n - p_n\|_{\alpha,\beta,p} \leq \left(\int_{-1}^{1} |w_{\alpha,\beta}(f_n - p_n)|^p \, dx\right)^{\frac{1}{p}}.$$

Theorem 2.11 immediately give

$$\mathcal{E}_n^{(2)}(f_n, w_{\alpha,\beta}, Y_s)_p \leq 2c(g,\mu)^{-1}\|f_n\|_{\alpha,\beta,p} + \omega_{i-1,i-1}^{\phi}\left(f_n^{(i-1)}, t\right)_{\alpha,\beta,p}. \tag{11}$$

Let

$$f_1(x) = 1, \quad f_2(x) = \frac{x - x^\#}{x_* - x_i} \text{ and } f_3(x) = \frac{x - x_*}{x_i - x^\#}, \tag{12}$$

therefore, (8), (11) and (12) implies

$$\mathcal{E}_n^{(2)}(f_n, w_{\alpha,\beta}, Y_s)_p \leq \mathcal{E}_n^{(2)}(f_1, w_{\alpha,\beta}, Y_s)_p + \mathcal{E}_n^{(2)}(f_2, w_{\alpha,\beta}, Y_s)_p + \mathcal{E}_n^{(2)}(f_3, w_{\alpha,\beta}, Y_s)_p$$

$$\leq 2c(g,\mu)^{-1}\left(\|f_1\|_{\alpha,\beta,p} + \|f_2\|_{\alpha,\beta,p} + \|f_3\|_{\alpha,\beta,p}\right) + \omega_{i-1,i-1}^{\phi}\left(f_1^{(i-1)}, t\right)_{\alpha,\beta,p}$$

$$+\omega_{i-1,i-1}^{\phi}\left(f_2^{(i-1)}, t\right)_{\alpha,\beta,p} + \omega_{i-1,i-1}^{\phi}\left(f_3^{(i-1)}, t\right)_{\alpha,\beta,p}$$

and

$$\mathbb{E}_n(f_n)_{\alpha,\beta,p} = \inf\left\{\|f_n - p_n\|_{\alpha,\beta,p}, \ p_n \in \pi_n \cap \mathcal{I}_f, \ f_n \in \mathbb{L}_p^{\alpha,\beta} \cap \Delta^{(2)}(Y_s) \cap \mathcal{I}_f\right\}$$

$$\leq \|f_n - p_n\|_{\alpha,\beta,p} \leq c(\alpha,\beta,p)\left(\int_{-1}^{1} |w_{\alpha,\beta}f_n|^p \, dx\right)^{\frac{1}{p}} \text{ where } c \text{ is large enough.}$$

From (12), therefore,

$$\mathbb{E}_n(f_n)_{\alpha,\beta,p} \leq c(\alpha,\beta,p)\left(\|f_1\|_{\alpha,\beta,p} + \|f_2\|_{\alpha,\beta,p} + \|f_3\|_{\alpha,\beta,p}\right).$$

Let us choose $\varepsilon > 0$, for any $\xi, \zeta > 0$, $\varepsilon < \xi$ and $\varepsilon < \zeta$. Then,

**Case I. CONVEXTIY:**

$$\mathcal{A} = \left\{n: \mathcal{E}_n^{(2)}(f_n, w_{\alpha,\beta}, Y_s)_p \geq \xi\right\}$$

$$\mathcal{A}_1 = \left\{n: \mathcal{E}_n^{(2)}(f_1, w_{\alpha,\beta}, Y_s)_p \geq \frac{\xi - \varepsilon}{3}\right\}$$

$$\mathcal{A}_2 = \left\{n: \mathcal{E}_n^{(2)}(f_2, w_{\alpha,\beta}, Y_s)_p \geq \frac{\xi - \varepsilon}{3}\right\}$$

$$\mathcal{A}_3 = \left\{n: \mathcal{E}_n^{(2)}(f_3, w_{\alpha,\beta}, Y_s)_p \geq \frac{\xi - \varepsilon}{3}\right\}.$$

Thus,

$$\mathcal{A} \subseteq \mathcal{A}_1 \cup \mathcal{A}_2 \cup \mathcal{A}_3,$$

whence

$$\sum_{n\in\mathcal{A}} a_n^{kj} \leq \sum_{n\in\mathcal{A}_1} a_n^{kj} \cup \sum_{n\in\mathcal{A}_2} a_n^{kj} \cup \sum_{n\in\mathcal{A}_3} a_n^{kj}.$$

**Case II. UNCONSTRAINED POLYNOMIAL:**

$$\mathcal{B} = \{n: \mathbb{E}_n(f_n)_{\alpha,\beta,p} \geq \zeta\}$$

$$\mathcal{B}_1 = \left\{n: \mathbb{E}_n(f_1)_{\alpha,\beta,p} \geq \frac{\zeta-\varepsilon}{3}\right\}$$

$$\mathcal{B}_2 = \left\{n: \mathbb{E}_n(f_2)_{\alpha,\beta,p} \geq \frac{\zeta-\varepsilon}{3}\right\}$$

$$\mathcal{B}_3 = \left\{n: \mathbb{E}_n(f_3)_{\alpha,\beta,p} \geq \frac{\zeta-\varepsilon}{3}\right\}.$$

Thus,

$$\mathcal{B} \subseteq \mathcal{B}_1 \cup \mathcal{B}_2 \cup \mathcal{B}_3,$$

whence

$$\sum_{n\in\mathcal{B}} b_n^{kj} \leq \sum_{n\in\mathcal{B}_1} b_n^{kj} \cup \sum_{n\in\mathcal{B}_2} b_n^{kj} \cup \sum_{n\in\mathcal{B}_3} b_n^{kj}.$$

By Cases I and II, (9) is proved. □

## 5. Data Availability

No data were used to support this study.

## 6. Conflict of Interests

The authors declare that there is no conflict of interests regarding the publication of this paper.

## 7. Acknowledgement

The first author is indebted to Prof. Eman Samir Bhaya (University of Babylon) for useful discussions of the subject. The first author is supported by University of Al-Muthanna while studying for his Ph. D. We would like to thank Universiti Utara Malaysia (UUM) for the financial support.

Malik Saad Al-Muhja

*Department of Mathematics and statistics, School of Quantitative Sciences, College of Arts and Sciences, Universiti Utara Malaysia, 06010 Sintok, Kedah, Malaysia.*
*Department of Mathematics and Computer Application, College of Sciences, University of Al-Muthanna, Samawa 66001, Iraq.*

E-mail addresses: dr.al-muhja@hotmail.com; malik@mu.edu.iq

Habibulla Akhadkulov

*Department of Mathematics and statistics, School of Quantitative Sciences, College of Arts and Sciences, Universiti Utara Malaysia, 06010 Sintok, Kedah, Malaysia.*

E-mail addresses: habibulla@uum.edu.my

Nazihah Ahmad

*Department of Mathematics and statistics, School of Quantitative Sciences, College of Arts and Sciences, Universiti Utara Malaysia, 06010 Sintok, Kedah, Malaysia.*

E-mail addresses: nazihah@uum.edu.my